\documentclass{elsart} 
\usepackage{amscd,amssymb,graphics}
\usepackage{mathrsfs} 
\usepackage[centertags]{amsmath}
\usepackage{latexsym}
\usepackage{amsfonts}
\usepackage{amssymb}

\newtheorem{theorem}{Theorem}[section]
\newtheorem{remark}{Remark}[section]

\newtheorem{lemma}{Lemma}[section]

\newtheorem{corollary}[theorem]{Corollary}

\newcommand{\abs}[1]{\vert#1\vert}
\def\norm#1{\left\Vert#1\right\Vert}
 
\def\R {{\mathbb R}}

\def\Ur{{\mathbb U}}
\def\Urm{{\mathbb U}_{\mathfrak m}}

\def\e{{\varepsilon}}

\def\Iso{{\mathrm{Iso}\,}}

\begin{document}
\begin{frontmatter}

\title{Subgroups of isometries of
Urysohn-Kat\v etov metric spaces of uncountable density}

\author{Brice R. Mbombo}
\address{Department of Mathematics, Faculty of Science,
University of Yaound\'e I,\\ PO Box 812, Yaound\'e, Cameroon}
\ead{bricero@yahoo.fr}
\author{Vladimir G. Pestov}
\address{Department of Mathematics and statistics, University of Ottawa,\\ 585 King Edward Avenue, Ottawa, Ontario, Canada K1N 6N5}
\ead{vpest283@uottawa.ca}

\begin{abstract}
According to Kat\v etov (1988),
for every infinite cardinal $\mathfrak m$ satisfying ${\mathfrak m}^{\mathfrak n}\leq {\mathfrak m}$ for all ${\mathfrak n}<{\mathfrak m}$, there exists a unique $\mathfrak m$-homogeneous universal metric space $\Ur_{\mathfrak m}$ of weight $\mathfrak m$. This object generalizes the classical Urysohn universal metric space $\Ur = \Ur_{\aleph_0}$. 
We show that for $\mathfrak m$ uncountable, the isometry group $\Iso(\Urm)$ with the topology of simple convergence is not a universal group of weight $\mathfrak m$: for instance, it does not contain $\Iso(\Ur)$ as a topological subgroup. More generally,
every topological subgroup of $\Iso(\Urm)$ having density $<{\mathfrak m}$ and possessing the bounded orbit property $(OB)$ is functionally balanced: right uniformly continuous bounded functions are left uniformly continuous. This stands in sharp contrast with Uspenskij's 1990 result about the group $\Iso(\Ur)$ being a universal Polish group.
\end{abstract}

\end{frontmatter}

\section{Introduction}
Some of the most interesting examples of non-locally compact topological groups studied recently are groups of automorphisms of countably infinite ultrahomogeneous Fra\"\i ss\'e structures \cite{KPT,KR} as well as their continuous analogues, including the group of isometries $\Iso(\Ur)$ of the Urysohn universal separable metric space $\Ur$ \cite{usp5,gro,NVT}.

There has been not much work done for non-separable versions of the above concepts. A notable exception is a classification by Lieberman \cite{lieberman} of unitary representations of the permutation group $S(\Gamma)$ of a set $\Gamma$ of an arbitrary infinite cardinality. Of course $S(\Gamma)$ is the automorphism group of a Fra\"\i ss\'e structure with empty signature. The aim of this note is to point out that for a non-trivial signature the uncountable case may pose new challenges.

Our case in point is the Urysohn universal metric space.
One of the central observations about this object is Uspenskij's result \cite{usp4} stating that the group $\Iso(\Ur)$ is a universal Polish group: every second-countable topological group is isomorphic with a suitable topological subgroup of $\Iso(\Ur)$. At the same time, the question of existence of a universal topological group of a given uncountable weight $\mathfrak m$ remains open to the day. In fact, it is open for {\em any} given cardinal $\mathfrak m>\aleph_0$.

In this connection, it is rather natural to begin by examining the group of isometries of a non-separable version of the Urysohn space constructed by Kat\v etov \cite{kat} under the following assupmtion on the cardinal $\mathfrak m$:
\begin{equation}
\label{eq:mn}
\sup\left\{{\mathfrak m}^{\mathfrak n}\colon {\mathfrak n}<{\mathfrak m}\right\}={\mathfrak m}.
\end{equation}
For instance, $\aleph_0$ and every strongly inaccessible cardinal $\mathfrak m$ are such.  
%
%
And GCH implies the above property for every infinite successor cardinal $\mathfrak m={\mathfrak n}^+$.

Assuming (\ref{eq:mn}), there exists a unique up to an isometry complete metric space $\Ur_{\mathfrak m}$ of weight $\mathfrak m$, which contains an isometric copy of every other metric space of weight $\leq\mathfrak m$ and is $<{\mathfrak m}$-homogeneous, that is, an isometry between any two metric subspaces of density $<{\mathfrak m}$ extends to a global self-isometry of $\Ur_{\mathfrak m}$. In particular, $\Ur_{\mathfrak \aleph_0}$ is just the classical Urysohn space $\Ur$ \cite{U25,U27}.

The topological group of isometries $\Iso(\Ur_{\mathfrak m})$, equipped with the topology of simple convergence, has weight $\mathfrak m$, and was suggested as a candidate for a universal topological group of weight $\mathfrak m$ by Iliadis (a private communication, 2006).

Somewhat surprisingly, it turns out not to be the case. Recall that a topological group $G$ has the property $(OB)$ (from the French for ``bounded orbits'') \cite{rosendal} if all the orbits of continuous actions of $G$ by isometries are bounded. We show that every subgroup $G$ of $\Iso(\Ur_{\mathfrak m})$ having density $<{\mathfrak m}$ and possessing the property $(OB)$ is necessarily a functionally balanced group: every left uniformly continuous bounded function on $G$ is right uniformly continuous. In particular, if such $G$ is either metrizable or locally connected, then its left and right uniformities coincide, or, equivalently, $G$ has a basis at identity consisting of conjugation-invariant open sets. (Topological groups with this property are known as {\em SIN groups.})

This is a serious restriction, which shows, in particular, that for uncountable $\mathfrak m$ the group $\Iso(\Ur_{\mathfrak m})$ does not even contain a copy of the topological group $\Iso(\Ur)$. In Uspenskij's terminology \cite{usp5}, an isometric embedding of the Urysohn space $\Ur$ into the space $\Ur_{\mathfrak m}$ (which is essentially unique up to a global isometry of $\Urm$) is not a $g$-embedding.

The techniques of our proof come from an earlier work \cite{PS,protasov,MNP} on another challenging open problem \cite{itz0,itz}: suppose every left uniformly continuous function on a topological group $G$ is right uniformly continuous, does it mean that $G$ is SIN? 

We prove that $\Iso(\Urm)$ contains every metrizable SIN group of weight $\leq\mathfrak m$, as well as every group of weight $\leq\mathfrak m$ in which the intersections of every family of $<\mathfrak m$ open subsets is open (a $P_{\mathfrak m}$-group). The problem of characterizing topological subgroups of the group $\Iso(\Urm)$ remains open.

\section{The construction of the space $\Urm$}

Here we will outline the main features of the space $\Urm$ as relevant to our results. The details of the construction can be found in \cite{kat} and \cite{usp5}. See also Chapter 3 in \cite{gro}, \cite{melleray}, \cite{NVT}, or Chapter V in \cite{P06} for the separable case.

A $1$-Lipschitz real-valued function $f$ on a metric space $X$ is a {\em Kat\v etov function} (although such functions were already thoroughly studied by Flood \cite{Fl1,Fl2}) if for every $x,y\in X$
\[d(x,y)\leq f(x)+f(y).\]
In particular, it follows that $f$ always takes non-negative values. It is easily seen that Kat\v etov functions are precisely the distance functions to points in metric extensions of $X$. In other words, $f\colon X\to\R$ is Kat\v etov if and only if there exist a metric space $Y$ containing $X$ as a subspace and a point $y\in Y$ with $f(x) = d(x,y)$ for each $x\in X$. (It may in particular happen that $y\in X$.)

Kat\v etov functions can always be extended from metric subspaces as follows. Let $X$ be a metric space and $A$ a metric subspace. Let $g\colon A\to\R$ be a Kat\v etov function. Then the formula
\begin{equation}
\label{eq:control}
g^{\uparrow X}(x) = \inf\{ g(a) + d(x,a)\colon a\in A\}\end{equation}
defines a Kat\v etov function $g^{\uparrow X}$ on $X$, known as the {\em Kat\v etov extension} of $g$ over $X$. This $g^{\uparrow X}$ is the largest among all $1$-Lipschitz functions on $X$ assuming prescribed values on $A$. For example, if $A=\{x_0\}$ is a singleton, then $\left(0_{x_0}\right)^{\uparrow X}=d(x_0,-)$ is the distance function to $x_0$.

If a Kat\^etov function $f$ on $X$ is of the form $f=  \left(f\vert_A\right)^{\uparrow X}$ for some $A\subseteq X$, then we will say that $f$ is {\em controlled} by $A$. In this case, the values of $f$ on $X$ can be restored by knowing the values of $f$ at the points of $A$.

For a metric space $X$, denote $E(X)$ the collection of all Kat\v etov functions on $X$ equipped with the supremum metric. Since for any $x,x_0\in X$ and $f,g\in E(X)$ one has $\abs{f(x)-g(x)}\leq f(x_0)+g(x_0)$, the metric is well-defined on $E(X)$ (that is, the supremum of $\abs{f(x)-g(x)}$ over $x\in X$ is always finite), unlike on the larger space of all $1$-Lipschitz functions. 

If $\mathfrak m$ is a cardinal number,  $E_{\mathfrak m}(X)$ denotes the subspace of $E(X)$ consisting of functions controlled by subspaces of cardinality (equivalently: density) $<{\mathfrak m}$. The density of $E_{\mathfrak m}(X)$ does not exceed $\sup\{d(X)^{\mathfrak n}\colon {\mathfrak n}<{\mathfrak m}\}$, where $d(X)$ denotes the density of $X$. 

The space $X$ isometrically embeds into $E_1(X)$ (and thus in every $E_{\mathfrak m}(X)$) in a canonical way, through the {\em Kuratowski embedding} $x\mapsto d(x,-)$. 

Suppose $X$ is a metric space of weight $\mathfrak m$. Consider the chain of iterated extensions $E_{\mathfrak m}^{\tau+1}(X) = E_{\mathfrak m}(E_{\mathfrak m}^{\tau}(X))$, where $\tau<\mathfrak m$ and $E_{\mathfrak m}^{\tau}(X) = \cup_{\lambda<\tau}E_{\mathfrak m}^{\lambda}(X)$ for limit cardinals.
Under the assumption (\ref{eq:mn}) on $\mathfrak m$, the density of every $E_{\mathfrak m}^{\tau}(X)$, $\tau<\mathfrak m$ is bounded by $\mathfrak m$. Thus,
iterating the construction $\mathfrak m$ times and taking the completion of the union
\begin{equation}
\label{eq:chain}
\bigcup_{{\mathfrak n}<{\mathfrak m}} E_{\mathfrak m}^{\mathfrak n}(X),\end{equation}
one obtains an embedding of $X$ into a complete metric space $\Ur_{\mathfrak m}$ of density $\mathfrak m$ which has the following version of the one-point extension property: if $Y$ is a subspace of $\Ur_{\mathfrak m}$ of density $<{\mathfrak m}$ and $f$ is a Katetov function on $Y$, then there exists a point $x\in \Urm$ with $f = d(x,-)$. A transfinite variation of a usual back-and-forth argument establishes that a space $\Urm$ with such properties is unique and $<{\mathfrak m}$-homogeneous. It follows from the construction that $\Ur_{\mathfrak m}$ is also universal for spaces of weight $\leq \mathfrak m$.
Notice that for $\mathfrak m\geq {\mathfrak c}$ (that is, under our assumption (\ref{eq:mn}), for uncountable $\mathfrak m$) the metric space union of the transfinite chain (\ref{eq:chain}) will be automatically complete. 

\section{Neighbourhoods of identity in subgroups of $\Iso(\Urm)$}

The group $\Iso(\Urm)$ of isometries of the space $\Urm$ is a topological group if equipped with the topology of simple convergence on $\Urm$. A neighbourhood basis at identity consists of all sets of the form
\[V[x_1,x_2,\ldots,x_n;\e]=\{g\in\Iso(\Urm)\colon d(x_i,gx_i)<\e,~i=1,2,\ldots,n\},\]
where $x_1,\ldots,x_n$ is a finite subset of $\Urm$ and $\e>0$. It is not difficult to verify that the weight of $\Iso(\Urm)$ is exactly $\mathfrak m$.

The following observation is crucial for our subsequent proof.

\begin{lemma}
\label{l:nbhds}
Let $\mathfrak m$ be an infinite cardinal satisfying (\ref{eq:mn}), and
let $G$ be a subgroup of $\Iso(\Urm)$ of density $<\mathfrak m$. The sets
\[V[x;\e]\cap G,~~x\in \Urm,~~\e>0\]
form a neighbourhood basis at the identity for the topology of $G$.
\end{lemma}

\begin{pf}
Without loss in generality replace $G$ with a dense subgroup of cardinality $<\mathfrak m$. 
Now let $x_1,x_2,\ldots,x_n\in  \Urm$ and $\e>0$ be arbitrary. We will choose $y\in\Urm$ and $\gamma>0$ so that $V[y;\gamma]\cap G\subseteq V[x_1,\ldots,x_n;\e]$. 

Let $D$ be the diameter of the set $X=\{x_1,x_2,\ldots,x_n\}$. Assume $D>0$, for otherwise there is nothing to be proved.
Choose $\gamma>0$ so that $\gamma\leq\e$, $\gamma<D/3$, and the open balls of radius $n\gamma$ around the points $x_i$, $i=1,2,\ldots,n$, are pairwise disjoint. The function $f(x_i)=D+i\gamma$, $i=1,2,\ldots,n$, is $1$-Lipschitz on $X$, and moreover a Kat\v etov function. Use the same letter $f$ to denote the Kat\v etov extension $\left(f\vert_X\right)^{\uparrow \Urm}$ over the entire space $\Urm$ (cf. Eq. (\ref{eq:control})).

Let $i=1,2,\ldots,n$. Suppose $x\in\Urm$ and for all $j=1,2,\ldots,i$ one has $d(x_j,x)\geq (i-j)\gamma$. Then for all $j=1,2,\ldots,n$ we have $d(x,x_j)+f(x_j)\geq D+i\gamma$ and so $f(x)\geq D+i\gamma$ according to Eq. (\ref{eq:control}). We obtain the following property of $f$: 
\begin{equation}
\label{eq:property}
\forall x\in\Urm,~~(f(x)<D+i\gamma)\Rightarrow \left(x\in \bigcup_{j=1}^{i-1} B_{(i-j)\gamma}(x_j)\right).
\end{equation}

Denote $A$ the union of $G$-orbits of the points $x_1,\ldots,x_n$. 
The one-point property of the space $\Urm$ implies the existence of a point $y\in\Urm$ with $d(y,a) =f(a)$ for each $a\in A$.
We claim that
\[V[y;\gamma]\subseteq V[x_1,\ldots,x_n;\e].\]
To see this, suppose $d(y,gy)<\gamma$. 
One has for every $x\in A$:
\[\abs{f(x)-f(gx)} = \abs{d(y,x) - d(y,gx)} = \abs{d(gy,gx)-d(y,gx)}
\leq d(y,gy) <\gamma.
\]
Applying this to $x_i$, $i=1,2,\ldots,n$, we obtain $f(gx_i)<D+(i+1)\gamma$ and consequently, modulo Eq. (\ref{eq:property}),
\[gx_i\in \bigcup_{j=1}^{i} B_{(i-j+1)\gamma}(x_j).\]
Since $g$ is an isometry, the choice of $\gamma$ assures that for $i\neq k$ the images $gx_i$ and $gx_k$ belong to two distinct balls as above. Assigning $j=\phi(i)$, where $j$ is the center of the ball $B_{(i-j+1)\gamma}(x_j)$ containing $gx_i$, we get an injection from $1,2,\ldots,n$ to itself with the property $\phi(i)\leq i$. This means $gx_i\in B_\e(x_i)$ for all $i$, as required.
\end{pf}

\section{Balanced and functionally balanced groups}

To fix the notation, we agree that the left uniform structure of a topological group $G$ has a basis of entourages $V_{L}=\{(x,y)\in G\times G:\,x^{-1}y\in V\}$ where $V$ runs over neighbourhoods of the identity. To obtain the right uniform structure, one replaces $x^{-1}y$ with $xy^{-1}$. The left and the right uniformities on a group $G$ coincide if and only if $G$ has a basis at identity consisting of open conjugation-invariant neighbourhoods (an easy check). Such topological groups are called {\em SIN groups} (from ``small invariant neighbourhoods''), or sometimes {\em balanced groups}. Examples of SIN groups include abelian and compact groups, as well as bounded groups of operators on Banach spaces equipped with the uniform operator topology.

A real-valued function $f$ on $G$ is said to be {\em left uniformly continuous} if it is uniformly continuous with regard to the left uniform structure. Equivalently, given $\e>0$, there is a neighbourhood of identity $V$ with $\abs{f(x)-f(y)}<\e$ whenever $x^{-1}y\in V$. A function $f$ is right uniformly continuous if the same conclusion holds as soon as $xy^{-1}\in V$. 

A topological group $G$ is called {\em functionally balanced,} or sometimes {\em FSIN} (``functionally SIN'') if every right uniformly continuous bounded function on $G$ is left uniformly continuous. 

Clearly, every balanced topological group is functionally balanced. The converse implication has been established for locally compact groups \cite{itz0}, metrizable groups \cite{protasov}, and locally connected groups \cite{MNP}, among others. Rather surprisingly, in the general case it remains an open problem, known as ``the Itzkowitz problem.'' For a recent survey with references, see \cite{BT}.

Here is a useful criterion of FSIN property.

\begin{theorem}[Protasov and Saryev \cite{PS}]\label{th:ps}
A topological group $G$ is functionally balanced if and only if for every subset $A\subseteq G$ and each neighbourhood of identity $V$ there is a neighbourhood of identity $U$ with $UA\subseteq AV$.
\end{theorem}

In the terminology of \cite{MNP}, this is the case where every subset $A$ of $G$ is a {\em neutral} set. We will somewhat strengthen the above result. A subset $A$ of a topological group $G$ is {\em left uniformly discrete} if it is uniformly discrete with regard to the left uniform structure. In other words, there is a neighbourhood of the identity $V$ such that, whenever $a,b\in A$ and $a\neq b$, one has $aV\cap bV=\emptyset$. Equivalently, there is a left-invariant continuous pseudometric $\rho$ on $G$ with regard to which the distances between pairs of distinct points in $A$ are uniformly bounded away from zero.

\begin{lemma}
A topological group $G$ is functionally balanced if and only if for every left uniformly discrete subset $A\subseteq G$ and each neighbourhood of identity $V$ there is a neighbourhood of identity $U$ with $UA\subseteq AV$.
\end{lemma}

\begin{pf}
It suffices to establish sufficiency of the condition, and we will accomplish this by appealing to Protasov--Saryev theorem \ref{th:ps}. 
Let $A$ be an arbitrary subset of $G$, and let $V$ be a neighbourhood of identity. Choose symmetric $W$ with $W^4\subseteq V$, and let $B$ be a maximal subset of $AW$ with the property that for all $a,b\in B$, $a\neq b$, the sets $aW$ and $bW$ are disjoint. Then clearly $A\subseteq BW^2$. Since $B$ is left uniformly discrete, by hypothesis there is a neighbourhood $U$ of the identity with the property $UB\subseteq BW$. Now one has
\[UA\subseteq UBW^2 \subseteq BW^3\subseteq AW^4\subseteq AV,\]
as required.
\end{pf}

\begin{remark}
\label{r:dense}
It is a trivial observation that if a functionally balanced topological group $G$ is dense in a group $H$, then $H$ is functionally balanced as well.
\end{remark}

We mention two of the known results about when functional balance implies the SIN property. 

\begin{theorem}[Protasov \cite{protasov}]
\label{th:protasov}
Every metrizable functionally balanced group is SIN.
\end{theorem}

\begin{theorem}[Megrelishvili--Nickolas--Pestov \cite{MNP}]
\label{th:mnp}
Every locally connected functionally balanced group is SIN.
\end{theorem}

For some further recent developments, see \cite{BT2,BB}.

\section{Functional balance of some groups of isometries of $\Urm$}

A topological group $G$ has {\em property} $(OB)$, or {\em topological Bergman property,} if every continuous left-invariant pseudometric on $G$ is bounded, or, equivalently, if every orbit of each continuous action of $G$ by isometries on a metric space $X$ is a bounded subset of $X$. Many examples of such groups can be found in \cite{rosendal}. They include, among others, all Bergman groups \cite{bergman}, such as the infinite permutation group $S_\infty$, even equipped with the discrete topology; the unitary group $U(\ell^2)$ with the uniform (hence also with the strong operator) topology; and the isometry group of the Urysohn sphere (that is, a sphere in the Urysohn space).

Here is our main result.

\begin{theorem}
\label{th:main}
Let $\mathfrak m$ be an uncountable cardinal satisfying the condition (\ref{eq:mn}). Denote by $\Iso(\Urm)$ the group of isometries of  the Urysohn--Kat\v etov metric space $\Urm$ of weight $\mathfrak m$ equipped with the topology of simple convergence on $\Urm$. Suppose $G$ is a topological subgroup of $\Iso(\Urm)$ of density $<\mathfrak m$, having property $(OB)$. Then $G$ is functionally balanced.
\end{theorem}

\begin{pf}
In view of Remark \ref{r:dense}, we can assume without loss in generality that the cardinality of $G$ does not exceed $\mathfrak m$.
Let $A\subseteq G$ be a left uniformly discrete subset of $G$, and let $V$ be an arbitrary neighbourhood of identity in $G$. We will find a neighbourhood $W$ of identity with the property $WA\subseteq AV$.

Without loss in generality, we can assume that $V$ is symmetric and that 
\begin{equation}
\label{eq:ab}
aV^2\cap bV^2=\emptyset
\end{equation}
whenever $a,b\in A$, $a\neq b$. 
Using Lemma \ref{l:nbhds}, select $x\in \Urm$ and $\e>0$ so that 
\[G\cap V[x;2\e]\subseteq V.\]
Let $(Ax)_{\e}$ denote the open $\e$-neighbourhood of the set $Ax$ in $\Urm$:
\[(Ax)_\e = \bigcup_{a\in A}B_\e(a).\]

Notice that if $a\neq b$, then open $\e$-balls around $ax$ and $bx$ do not intersect. Indeed, if $B_\e(ax)\cap B_\e(bx)\neq \emptyset$, then 
\[d(b^{-1}ax,x)=d(ax,bx)<2\e,\]
whence $b^{-1}a\in V[x;2\e]\cap G\subseteq V$ and $a\in bV$. This implies $a=b$ in view of (\ref{eq:ab}). 

Let $D$ be a positive number greater than both the diameter of $Ax$ and $4\e$.
Let $F$ stand for the complement to $(Ax)_{\e}$ in $\Urm$. Now
define a $1$-Lipschitz function $f\colon\Urm\to \R$ as follows:
\[f(y) = D - d(F,y).\]
Since the open $\e$-balls around elements of $A$ are pairwise disjoint, the values of $f$ are bounded from below by $D-2\e$. Since $2(D-2\e)> D$, it follows that $f$ is a Kat\^etov function. The one-point extension property of $\Urm$ modulo the fact that the $G$-orbit of $x$ has cardinality $\leq\mathfrak m$ assures the existence of a point $z\in\Urm$ such that
\[\forall g\in G,~~f(gx)=d(gx,z).\]
Now set
\[W = V[z;\e/3]\cap G.\]

Let $w\in W$ and $a\in A$ be arbitrary. Since $w$ only moves $z$ by a distance strictly less than $\e/3$, the triangle inequality assures that the values of $f$ on $wax$ and on $ax$ differ by strictly less than $\e$. One has $f(ax)=D-\e$, and therefore $f(wax)<D$. We conclude: $wax\in (Ax)_{\e}$. Put otherwise, there is $b\in A$ with $d(wax,bx)<\e$, or equivalently, $wax\in B_\e(bx)$. This in turn means that $wa\in AV[x;\e]\subseteq AV$, as required.
\end{pf}

Theorems \ref{th:protasov} of Protasov and \ref{th:mnp} of Megrelishvili--Nickolas--Pestov imply:

\begin{corollary}
Let $\mathfrak m$ be an uncountable cardinal satisfying (\ref{eq:mn}), and let $G$ be a topological subgroup of $\Iso(\Urm)$ of density $<{\mathfrak m}$ having property $(OB)$ which is either metrizable or locally connected. Then $G$ is a SIN group.
\end{corollary}

Now we obtain examples of Polish topological groups which admit no embedding into the group $\Iso(\Urm)$ for $\mathfrak m$ uncountable.

1. Denote $\Ur^{\bigcirc}$ the unit sphere in the Urysohn metric space.
The group $\Iso(\Ur^{\bigcirc})$ is both metrizable and locally connected (in fact, homeomorphic to $\ell^2$ \cite{melleray2}), has property $(OB)$ (\cite{rosendal}, Theorem 1.5) and is not SIN, which of course follows from the fact that it is a universal Polish group \cite{usp4,usp5}, but is most easily verified directly. Thus, $\Iso(\Ur^{\bigcirc})$ is not contained in $\Iso(\Urm)$.

2. The group $S_\infty$ of all permutations of a countably infinite set $\omega$, equipped with the topology of simple convergence on $\omega$ with a discrete topology. This group is Polish and has Bergman's property, hence the property $(OB)$. 

3. The group ${\mathrm{Homeo}}([0,1]^{\aleph_0})$ of homeomorphisms of the Hilbert cube equipped with the topology of compact convergence, because it is a non-SIN Polish group with property $(OB)$ \cite{rosendal}.

4. The group $U(\ell^2)$ with the strong operator topology has property $(OB)$ \cite{atkin}, is Polish and is not SIN.

On the contrary, note that the same group $U(\ell^2)$ with the {\em uniform} operator topology is a metrizable SIN group of weight continuum, and so embeds into $\Iso(\Urm)$ for $\mathfrak m$ uncountable by force of our Theorem \ref{th:sin} below.

5. It follows in particular that the (universal Polish) group $\Iso(\Ur)$ of isometries of the Urysohn metric space admits no embedding into $\Urm$ as a topological subgroup.

\section{Concluding remarks}

We believe that this work just lifts a small corner of the curtain that conceals an unusual nature of the topological automorphisms groups of ultrahomogeneous structures in the non-separable case. 

It would be interesting to characterize topological subgroups of the group $\Iso(\Urm)$ for uncountable $\mathfrak m$. Here are two additional results in this direction. 

\begin{theorem}
\label{th:sin}
Let $\mathfrak m$ to be an infinite cardinal satisfying (\ref{eq:mn}). Every metrizable SIN group of weight $\leq \mathfrak m$ embeds into $\Iso(\Urm)$.
\end{theorem}

The proof employs Kat\v etov functions in the same way as in Uspenskij's classical argument for embeddings into $\Iso(\Ur)$ \cite{usp4,usp5}. Namely, the action of $G$ is recursively extended from $E_{\mathfrak m}^{\tau}(X)$ to $E_{\mathfrak m}^{\tau+1}(X)$, and the continuity of the extension is assured by the following  observation.

\begin{lemma}
Let $G$ be a group acting by isometries on a metric space $X$, let $V$ be a neighbourhood of identity in $G$ and $\e>0$. Suppose the following property holds: 

\begin{center}
$(\star)$ whenever $x\in X$ and $v\in V$, one has $d(x,vx)\leq\e$. 
\end{center}

Then $(\star)$ also holds for the action of $G$ on $E(X)$ (the left regular representation defined by $^gf(x)=f(g^{-1}x)$).
\end{lemma}

\begin{pf}
Let $f$ be a $1$-Lipschitz function on $X$. One has 
\begin{eqnarray*}
\norm{^gf- f} &=& \sup_{x\in X} \abs{f(g^{-1}x)-f(x)} \\
&\leq&\sup_{x\in X} d(g^{-1}x,x) \\
&\leq & \e,
\end{eqnarray*}
provided $g\in V$.
\end{pf}

To finish the proof of Theorem \ref{th:sin}, fix a bi-invariant metric on a metrizable SIN group $G$ which generates the topology. For every value of $\e>0$, the property $(\star)$ holds for the action of $G$ on itself by left multiplication, with $V=B_\e(e)$. The same is therefore true for the iterated action of $G$ on $\Urm=\cup_{{\mathfrak n}<{\mathfrak m}}E_{\mathfrak m}^{({\mathfrak n})}(X)$, and it follows that this action determines a topological group embedding of $G$ into $\Iso(\Urm)$.

We do not know whether every SIN group of weight $\leq {\mathfrak m}$ embeds into $\Iso(\Urm)$. At the same time, not all subgroups of $\Iso(\Urm)$ are SIN, as follows from our next result.

\begin{theorem}
\label{th:emb}
Let $\mathfrak m$ to be an infinite cardinal satisfying (\ref{eq:mn}). Suppose a topological group $G$ of weight $\mathfrak m$ has the property that the intersection of every family of strictly less than $\mathfrak m$ open subsets is open. (Sometimes such groups are called $P_{\mathfrak m}$-groups.)
Then $G$ is isomorphic with a topological subgroup of $\Iso(\Urm)$. 
\end{theorem}
  
The step of the recursion is provided by the following lemma.

\begin{lemma}
Let $G$ be a topological group satisfying the assumptions of Theorem \ref{th:emb}. Suppose $G$ acts continuously on a metric space $X$ by isometries. Then the canonical extension of the action to $E_{\mathfrak m}(X)$ is continuous.
\end{lemma}

\begin{pf}
Let $f$ be a $1$-Lipschitz function on $X$ controlled by a subset $A$ of cardinality $<{\mathfrak m}$. Given an $\e>0$, set 
\[V=\bigcap_{a\in A}V[a;\e].\]
This $V$ is a neighbourhood of identity in $G$ by force of our assumption on the group, and if $g\in V$, then for each $a\in A$ one has $\abs{f(ga)-f(a)}<\e$. This condition, in its turn, assures that for every $x\in X$ and $g\in V$ one has $\abs{f(gx)-f(x)}<\e$, and so the action of $G$ on $E_{\mathfrak m}(X)$ is continuous.
\end{pf}

Non-discrete \cite{CR} and even non-SIN \cite{vp1} topological groups as above exist and have been used to produce counter-examples with various exotic combinations of properties.

One may ask: is the group $\Iso(\Urm)$ (path) connected, just like the group $\Iso(\Ur)$ is (cf. \cite{melleray2})? Locally connected?

The group $\Iso(\Urm)$ is extremely amenable (has the fixed point on compacta property), as follows from results of \cite{pestov02} (Theorems 6.5 and 6.6). It would be interesting to verify whether the space $\Urm$ has the following, closely related, property. An ultrahomogeneous (that is, $<\omega$-homogeneous) metric space $X$ is {\em oscillation stable} (see Chapter 8 in \cite{P06}) if every finite partition $\gamma$ of $X$ contains an element $A\in\gamma$ whose each $\e$-neighbourhood $A_\e$, $\e>0$, contains an isometric copy of $X$. A famous result of Odell and Schlumprecht \cite{OS} says that the unit sphere in the Hilbert space is not oscillation stable (has distortion), while a recent remarkable theorem of Nguyen Van Th\'e and Sauer \cite{NVTS} (based on an earlier work of Lopez-Abad and Nguyen Van Th\'e \cite{LANVT}) establishes the oscillation stability of the unit sphere of the separable Urysohn space $\Ur$. Is the unit sphere of the space $\Urm$, $\mathfrak m>\aleph_0$ oscillation stable?

Further, Hjorth has proved \cite{hjorth} that no Polish group is oscillation stable when it is acting on itself by left translations, and his proof was considerably simplified by Melleray \cite{melleray3} (for relevant definitions, see either any of the above papers, or else Chapter 8 in \cite{P06}). It remains unknown whether the result holds true for arbitrary topological groups, and the group $\Iso(\Urm^{\bigcirc})$ ($\mathfrak m$ uncountable) looks like a natural candidate to explore in this regard. (Here, as before, $\Urm^{\bigcirc}$ denotes the unit sphere in $\Urm$.)

One can also study non-separable analogues of the Rado universal graph $R$ (that is, a universal ultrahomogeneous metric space whose metric takes values in $\{0,1,2\}$), as well as other types of Fra\"\i ss\'e structures. Their groups of automorphisms can serve, in our view, as potential candidates for an example of a functionally balanced group that is not SIN.

One way to obtain non-separable ultrahomogeneous Fra\"\i ss\'e structures is through the ultraproduct construction. For instance, under CH the metric ultrapower of the Urysohn space $\Ur$ with regard to any non-principal ultrafilter on the natural numbers (cf. \cite{BYBHU}) is easily shown to be isometric to the space $\Ur_{\mathfrak c}$. (For closely related, deeper results, see \cite{FS}). When do groups of automorphisms of saturated non-separable continuous logic models of various metric structures exhibit different behaviour from the separable case? One important example where this does not happen is that of the unit sphere of a separable Hilbert space.

\section*{Acknowedgements}

This project was supported by the Natural Sciences and Engineering Research Council of Canada through the 2007--2012 NSERC Discovery grant ``Geometry of infinite-dimensional groups'', as well as by the University of Ottawa through internal grants. The authors are grateful to Ahmed Bouziad and Michael Megrelishvili for careful reading of earlier versions of the article and bringing our attention to some inaccuracies.


\begin{thebibliography}{100}

\bibitem{atkin} 
C.J. Atkin, {\em Boundedness in uniform spaces, topological groups, and homogeneous spaces,} Acta Math. Hungar. \textbf{57} (1991), 213--232.

\bibitem{BB} A. Bareche and A. Bouziad, 
{\em Itzkowitz's problem for groups of finite exponent,}
Topology Proc. \textbf{36} (2010), 375--383. 

\bibitem{BYBHU} I. Ben Yaacov, A. Berenstein, C.W. Henson, and A. Usvyatsov, {\em Model Theory for Metric Structures,} in: Model theory with applications to algebra and analysis. Vol. 2, 315--427, 
London Math. Soc. Lecture Note Ser., \textbf{350}, Cambridge Univ. Press, Cambridge, 2008.

\bibitem{bergman}
G.M. Bergman, {\em Generating infinite symmetric groups,}
Bull. London Math. Soc. \textbf{38} (2006), 429--440. 

\bibitem{BT2} A. Bouziad and J.-P. Troallic, {\em
Left and right uniform structures on functionally balanced groups,}
Topology Appl. \textbf{153} (2006), 2351--2361. 

\bibitem{BT} A. Bouziad and J.-P. Troallic, {\em Problems about the uniform structures of topological groups,} in: Open problems in topology. II,
Edited by Elliott Pearl, Elsevier B.V., Amsterdam, 2007, pp. 359--366.

\bibitem{CR}
W.W. Comfort and K.A. Ross, {\em Pseudocompactness and uniform continuity in topological groups,} Pacific J. Math. \textbf{16} (1966), 483--496.

\bibitem{FS} I. Farah and S. Shelah, {\em A dichotomy for the number of ultrapowers,} arXiv:0912.0406v1 [math.LO].

\bibitem{Fl1} J. Flood,
{\it Free Topological Vector Spaces,} Ph.D. thesis, 
Australian National University, Canberra, 1975, 109 pp.

\bibitem{Fl2} J. Flood,
{\it Free locally convex spaces,} Dissertationes Math.
{\bf CCXXI} (1984), PWN, Warczawa.

\bibitem{gro} M. Gromov, {\it Metric Structures for Riemannian and Non-Riemannian Spaces,} Progress in Mathematics 152, Birkhauser Verlag, 1999.

\bibitem{hjorth}
G. Hjorth, {\em An oscillation theorem for groups of isometries,} Geom. Funct. Anal. \textbf{18} (2008), 489--521.

\bibitem{itz0}
G. Itzkowitz, {\em Continuous measures, Baire category, and uniform continuity in topological groups,} Pacific J. Math. \textbf{54} (1974), 115--125.

\bibitem{itz} G. Itzkowitz, {\it Uniformities and uniform continuity on topological groups,} General topology and its Applications, {\bf 134}
(1991), 155--178.

\bibitem{itz2} G. Itzkowitz, {\it Projective limits and balanced topological groups,} General topology and its Applications, {\bf 110}
(2001), 163--183.

\bibitem{kat} M. Kat\v etov, {\it On universal metric spaces,} in: Gen. Topology and its Relations to Modern Analysis and Algebra VI, Proc. Sixth Prague Topol. Symp. 1986, Z. Frolik, ed., Heldermann Verlag (1988), 323--330.

\bibitem{KPT} A.S. Kechris, V.G. Pestov and S. Todorcevic,
\textit{Fra\"\i ss\'e limits, Ramsey theory, and topological dynamics of
automorphism groups,} Geom. Funct. Anal. \textbf{15} (2005), 
106--189.

\bibitem{KR} A.S. Kechris and C. Rosendal, 
\textit{Turbulence, amalgamation and generic automorphisms of 
homogeneous structures,} Proc. Lond. Math. Soc. (3) \textbf{94} (2007), 302--350. 

\bibitem{lieberman}
A. Lieberman, {\em The structure of certain unitary representations of infinite symmetric groups,} Trans. Amer. Math. Soc. \textbf{164} (1972), 189--198. 

\bibitem{LANVT}
J. Lopez-Abad and L. Nguyen Van Th\'e, {\em The oscillation stability problem for the Urysohn sphere: a combinatorial approach,} Topology Appl. \textbf{155} (2008), 1516--1530.

\bibitem{MNP} M. Megreslishvili, P. Nickolas and V. Pestov, {\it Uniformities and uniformly continuous functions on locally connected groups,} Bull. Austral. Math. Soc. {\bf 56} (1997), 279--283.

\bibitem{melleray}
J. Melleray, {\em On the geometry of Urysohn's universal metric space,} Topology Appl. \textbf{154} (2007), 384--403.

\bibitem{melleray2}
J. Melleray, {\em Topology of the isometry group of the Urysohn space,} Fund. Math. \textbf{207} (2010), 273--287.

\bibitem{melleray3} J. Melleray, {\em A note on Hjorth's oscillation theorem,} Journal of Symbolic Logic, to appear, {\tt http://math.univ-lyon1.fr/$\sim$melleray/Oscillation.pdf} (accessed on Dec. 30, 2010).

\bibitem{OS}
E. Odell and T. Schlumprecht, \textit{The distortion problem,}
Acta Math. {\bf 173} (1994), 259--281.

\bibitem{NVTS}
L. Nguyen Van Th\'e and N.W. Sauer, {\em The Urysohn sphere is oscillation stable,} Geom. Funct. Anal. \textbf{19} (2009), 536--557.

\bibitem{NVT} L. Nguyen Van Th\'e, \textit{Structural Ramsey Theory of Metric Spaces and Topological Dynamics of Isometry Groups}, Mem. Amer. Math. Soc. \textbf{206} (2010), no. 968.

\bibitem{vp1} V.G. Pestov, {\em Embeddings and condensations of topological groups} (Russian), Mat. Zametki \textbf{31} (1982), no. 3, 443--446.

\bibitem{pestov02} V.G. Pestov, {\em Ramsey-Milman phenomenon, Urysohn metric spaces, and extremely amenable groups,} Israel Journal of Mathematics \textbf{127} (2002), 317--358. Corrigendum, ibid., \textbf{145} (2005), 375--379.

\bibitem{P06} V. Pestov, {\em Dynamics of Infinite-Dimensional Groups: the Ramsey-Dvoretzky-Milman phenomenon,} Amer. Math. Soc. University Lecture Series \textbf{40}, 2006.

\bibitem{protasov}
I.V. Protasov, {\em Functionally balanced groups,}
Math. Notes \textbf{49} (1991), 614--616. 

\bibitem{PS} I.V. Protasov and A. Saryev, {\em The semigroup of closed subsets of a topological group} (Russian),
Izv. Akad. Nauk Turkmen. SSR Ser. Fiz.-Tekhn. Khim. Geol. Nauk 1988, no. 3, 21--25. 

\bibitem{rosendal} C. Rosendal, {\em A topological version of the Bergman property,} Forum Math. \textbf{21} (2009), 299--332. 

\bibitem{U25}
P. S. Urysohn, \textit{Sur un espace m\'etrique universel,}
C. R. Acad. Sci. Paris \textbf{180} (1925), 803--806.

\bibitem{U27} 
P. Urysohn, \textit{Sur un espace m\'etrique universel,} Bull. Sci.
Math. {\bf 51} (1927), 43--64 et 74--90.
\bibitem{usp4} V.V. Uspenskij, {\it On the group of isometries of the Urysohn universal metric space,} Comment. Math. Univ. Carolinae {\bf 31} (1990), 181--182.
 
\bibitem{usp5} V.V. Uspenskij, 
{\em On subgroups of minimal topological groups,}
Topology Appl. \textbf{155} (2008), 1580--1606.

\end{thebibliography}
\end{document}